\documentstyle[12pt,leqno,amstex,amssymb,amsfonts,eucal]{article}
\begin{document}
\oddsidemargin= 5mm
\topmargin= 35mm
\textwidth=170mm
\textheight=270mm
\pagestyle{plain}
\newcounter{r}
\newcommand{\Ker}{Ker}
\newcommand{\Gal}{Gal}
\newcommand{\Mat}{Mat}
\newcommand{\sign}{sign}
\newcommand {\Log}{Log}
\newcommand{\Res}{Res}
\newcommand{\Tr}{Tr}
\newcommand{\Nm}{Nm}
\begin{center}
{\large\bf  ON SOME SYSTEMS OF DIFFERENCE EQUATIONS }
\end{center}
\begin{center}
{\large\bf L.A.Gutnik}
\end{center}
\vskip.10pt

 {\hfill\sl Dedicated to the memory of}

\hspace*{8cm}{\hfill\sl  Professor N.M.Korobov.}
\vskip.15pt
Let
$$
\vert z\vert\ge1,-3\pi/2<\arg(z)\le\pi/2,\log(z)=\ln(\vert z\vert)+i\arg(z).
$$
Then $\log(-z)=\log(z)-i\pi,$ if $\Re(z)>0$
 and $\log(z)=\log(-z)-i\pi,$ if $\Re(z)<0.$

Let
\begin{equation}\label{eq:b}
f_{l,1}(z, \nu) =
\sum\limits_{k=0}^\nu(-1)^{(\nu+k)l}(z)^k\binom{\nu}k^{2+l}
\binom{\nu+k}\nu^{2+l},
\end{equation}
where $l=0,\,1,\,2,\,\nu\in[0,+\infty)\cap{\mathbb Z}.$
Let
\begin{equation}\label{eq:d}
R(t,\nu)=\frac{\prod\limits_{j=1}^\nu(t-j)}
{\prod\limits_{j=0 }^\nu(t+j)}. 
\end{equation}
where $\nu\in[0,+\infty)\cap{\mathbb Z},$
\begin{equation}\label{eq:g}
f_{l,2}(z, \nu)=\sum\limits_{t=1+\nu}^{+\infty}z^{-t }(R(t,\nu))^{2+l},
\end{equation}
where $l=0,\,1,\,2$ and $\nu\in[0,+\infty)\cap{\mathbb Z},$ and
since $(R(t,\nu))^{2+l}$ for $\nu\in\Bbb N$
has in the points $t = 1,\,\ldots,\,\nu,$ the zeros of the order $2+l,$
it follows that
\begin{equation}\label{eq:h}
f_{l,2}(z,\nu)=\sum\limits_{t=1}^{+\infty}z^{-t}(R(t,\nu))^{2+l},
\end{equation}
for $l=0,\,1,\,2$ and $\nu\in[0,+\infty)\cap{\mathbb Z}.$
Let
\begin{equation}\label{eq:ac}
f_{l,3}(z,\nu)=(\log(z))f_{l,2}(z,\nu)+f_{l,4}(z,\nu),
\end{equation}
where
\newpage
\pagestyle{headings}
\topmargin= -15mm
\textheight=250mm
\markright{\footnotesize\bf L.A.Gutnik, On some systems of difference 
equations.}
\begin{equation}\label{eq:ab}
f_{l,4}(z,\nu)=-\sum\limits_{t=1+\nu}^{+\infty} z^{-t}
\left (\frac\partial{\partial t}(R^{2+l})\right)(t,\nu), 
\end{equation}
$l=0,\,1,\,2$ and $\nu\in[0,+\infty)\cap{\mathbb Z},$
 and since $(R(t,\nu))^{2+l}$ for $\nu\in\Bbb N$
has in the points $t = 1,\,\ldots,\,\nu,$ the zeros of the order $2+l,$
it follows that
\begin{equation}\label{eq:ab1}
f_{l,4}(z,\nu)=-\sum\limits_{t=1}^{+\infty} z^{-t}
\left (\frac\partial{\partial t}(R^{2+l})\right)(t,\nu)
\end{equation}
for $l=0,\,1,\,2$ and $\nu\in[0,+\infty)\cap{\mathbb Z}.$
Let
\begin{equation}\label{eq:ah}
f_{l,5}^\vee(z,\nu)=-i\pi f_{l,3}(z, \nu)+f_{l,5}(z, \nu),\end{equation}
with
$l=1,\,2,\,\nu\in[0,+\infty)\cap{\mathbb Z}$ and
\begin{equation}\label{eq:ai}
f_{l,5}(z,\nu)=
\end{equation}
$$2^{-1}(\log(z))^2f_{l,2}(z,\nu)+(\log(z))f_{l,4}(z,\nu)+f_{l,6}(z,\nu)=$$
$$=-2^{-1}(\log(z))^2f_{l,2}(z,\nu)+
+(\log(z))f_{l,3}(z,\nu)+
f_{l,6}(z,\nu),$$
where
\begin{equation}\label{eq:ag1}
f_{l,6}(z,\nu)=2^{-1}
\sum\limits_{t=1+\nu}^\infty z^{-t}
\left(\left(\frac\partial{\partial t}\right)^2(R^{2+l})\right)(t,\nu),
\end{equation}
 and
since $(R(t,\nu))^{2+l}$ for $\nu\in\Bbb N$
has in the points $t = 1,\,\ldots,\,\nu,$ the zeros of the order $2+l,$
and $l=1,\,2$ now,
it follows that
\begin{equation}\label{eq:ag2}
f_{l,6}(z,\nu)=2^{-1}
\sum\limits_{t=1+\nu}^\infty z^{-t}
\left(\left(\frac\partial{\partial t}\right)^2(R^{2+l})\right)(t,\nu)
\end{equation}
for $l=1,\,2$ and $\nu\in[0,+\infty)\cap{\mathbb Z}.$
Let
\begin{equation}\label{eq:bd}
f_{l,7}^\vee(z,\nu)=f_{l,7}(z,\nu)+(2\pi^2/3)f_{l,3}(z, \nu).
\end{equation}
with
$l=2,\,\nu\in[0,+\infty)\cap{\mathbb Z}$ and
\begin{equation}\label{eq:bb}
f_{l,7}(z,\nu)=
\end{equation}
$$-3^{-1}(\log(z))^3f_{l,2}(z,\nu)+2^{-1}(\log(z))^2f_{l,3}(z,\nu)+
f_{l,8}(z,\nu)+$$
$$(\log(z))(f_{l,5}(z,\nu)+
2^{-1}(\log(z))^2f_{l,2}(z,\nu)-(\log(z))f_{l,3}(z,\nu))=$$
$$6^{-1}(\log(z))^3f_{l,2}(z,\nu)-2^{-1}(\log(z))^2f_{l,3}(z,\nu)
+(\log(z))f_{l,5}(z,\nu)+f_{l,8}(z,\nu),$$
where
\begin{equation}\label{eq:bc}
f_{l,8}(z,\nu)=-6^{-1}
\sum\limits_{t=\nu+1}^\infty z^{-t}
\left(\left(\frac\partial{\partial t}\right)^3(R^{2+l})\right)(t,\nu),
\end{equation}
 and, since $(R(t,\nu))^{2+l}$ for $\nu\in\Bbb N$  
have in the points $t = 1,\,\ldots,\,\nu,$ the zeros of the order $2+l,$ 
 and $l=2$ now, it follows that
\begin{equation}\label{eq:bc1}
f_{l,8}(z,\nu)=-6^{-1}
\sum\limits_{t=1}^\infty z^{-t}
\left(\left(\frac\partial{\partial t}\right)^3(R^{2+l})\right)(t,\nu).
\end{equation}
Let
$${\mathfrak K}_0=\{1,\,2,\,3\},\,{\mathfrak K}_1=\{1,\,2,\,3,\,5\},\,
{\mathfrak K}_2=\{1,\,2,\,3,\,5,\,7\}.$$
Let $\lambda$ be a variable. We denote by $T_{n,\lambda}$ the
diagonal $n\times n$-matrix, $i$-th diagonal element of which
is equal to $\lambda^{i-1}$ for $i=1,\,...,\,n.$
We denote by $\delta$ the operator $z\frac{d}{dz}.$ Let further
$l=0,\,1,\,2,\,k\in{\mathfrak K}_l,\,\vert z\vert>1,\nu\in{\mathbb N},$
and let $Y_{l,k}(z;\nu)$ be the columnn with $4+2l$ elements, $i$-th of which
is equal to $(\nu^{-1}\delta)^{i-1}f^\vee_{l,k}(z,\nu)$ for $i=1,\,...,\,4+2l.$

\bf Theorem 1. \it The following equalities hold
\begin{equation}\label{eq:4d0}
A_l^\sim(z;\nu)Y_{l,k}(z;\nu)=T_{4+2l,1-\nu^{-1}}Y_{l,k}(z;\nu-1),
\end{equation}
\begin{equation}\label{eq:4d1}
Y_{l,k}(z;\nu)=
T_{4+2l,-1}A_l^\sim(z;-\nu)T_{4+2l,-1+\nu^{-1}}Y_{l,k}(z;\nu-1),
\end{equation}
where
$l=0,\,1,\,2,\,k\in{\mathfrak K}_l,\,\vert z\vert>1,\nu\in{\mathbb N},\nu\ge2,$
$$A_l^\sim(z;\nu)=S^\sim_l+z\sum\limits_{i=0}^{1+l}\nu^{-i}V_l^{\sim\ast}(i)$$
with
\begin{equation}\label{eq:4ab}
S^\sim_0=
\left(\matrix
1& {-4} &8&{-12}\\
0& {1} &{-4}&8\\
0& {0} &{1}&{-4}\\
0&  0  & 0  &1
\endmatrix\right)
\end{equation}
\begin{equation}\label{eq:4ac}
S^\sim_1=
\left(\matrix
{-1}&6&{-18}&{38}&{-66}&{102}\\
0&{-1}&6&{-18}&{38}&{-66}\\
0&0&{-1}&6&{-18}&{38}\\
0&0&{0}&{-1}&6&{-18}\\
0&0&0&0&{-1}&6&\\
0&0&0&0&0&{-1}
\endmatrix\right),
\end{equation}
\begin{equation}\label{eq:4ad}
S^\sim_2=
\left(\matrix
1&{-8}&{32}&{-88}&{192}&{-360}&{608}&{-952}\\
0&1&{-8}&{32}&{-88}&{192}&{-360}&{608}\\
0&0&1&{-8}&{32}&{-88}&{192}&{-360}\\
0&0&{0}&1&{-8}&{32}&{-88}&{192}\\
0&0&0&0&1&{-8}&{32}&{-88}\\
0&0&0&0&0&1&{-8}&{32}\\
0&0&0&0&0&0&1&{-8}\\
0&0&0&0&0&0&0&1
\endmatrix\right),
\end{equation}
$$V^{\sim\ast}_0(0)=4\left(\matrix
{4}&{-5}&{-2}&{3}\\
{-3}&{4}&1&{-2}\\
2&{-3}&0&1\\
{-1}&2&{-1}&0
\endmatrix\right),$$
$$V^{\sim\ast}_0(1)=4\left(\matrix
{3}&{-6}&{3}&{0}\\
{-2}&{4}&{-2}&{0}\\
1&{-2}&1&0\\
{0}&0&{0}&0
\endmatrix\right),
$$
$$V_1^{\sim\ast}(0)=\left(\matrix
{146}&{-198}&{-180}&{268}&{66}&{-102}\\
{-102}&{146}&{108}&{-180}&{-38}&{66}\\
{66}&{-102}&{-52}&{108}&{18}&{-38}\\
{-38}&{66}&{12}&{-52}&{-6}&{18}\\
{18}&{-38}&{12}&{12}&{2}&{-6}\\
{-6}&{18}&{-20}&{12}&{-6}&{2}
\endmatrix\right),$$
$$V_1^{\sim\ast}(1)=\left(\matrix
{240}&{-516}&{108}&{372}&{-204}&0\\
{-160}&{348}&{-84}&{-236}&{132}&0\\
{96}&{-212}&{60}&{132}&{-76}&0\\
{-48}&{108}&{-36}&{-60}&{36}&0\\
{16}&{-36}&{12}&{20}&{-12}&0\\
{0}&{-4}&{12}&{-12}&{4}&0
\endmatrix\right),$$
$$V_1^{\sim\ast}(2)=\left(\matrix
{102}&{-306}&{306}&{-102}&0&0\\
{-66}&{198}&{-198}&{66}&0&0\\
{38}&{-114}&{114}&{-38}&0&0\\
{-18}&{54}&{-54}&{18}&0&0\\
{6}&{-18}&{18}&{-6}&0&0\\
{-2}&{6}&{-6}&{2}&0&0
\endmatrix\right),$$
$$V_2^{\sim\ast}(0)=
8\left(\matrix
{176}&{-249}&{-364}&{545}&{280}&{-431}&{-76}&{119}\\
{-119}&{176}&{227}&{-364}&{-169}&{280}&{45}&{-76}\\
{76}&{-119}&{-128}&{227}&{92}&{-169}&{-24}&{45}\\
{-45}&{76}&{61}&{-128}&{-43}&{92}&{11}&{-24}\\
{24}&{-45}&{-20}&{61}&{16}&{-43}&{-4}&{11}\\
{-11}&{24}&{-1}&{-20}&{-5}&{16} &1&{-4}\\
4&{-11}&8&{-1}&4&{-5}&0&1\\
{-1}&4&{-7}&8&{-7}&4&{-1}&0
\endmatrix\right),$$
$$V_2^{\sim\ast}(1)=8\left(\matrix
{455}&{-1020}&{-113}&{1552}&{-603}&{-628}&{357}&0\\
{-300}&{682}&{44}&{-996}&{404}&{394}&{-228}&0\\
{185}&{-428}&{-3}&{592}&{-253}&{-228}&{135}&0\\
{-104}&{246}&{-16}&{-316}&{144}&{118}&{-72}&0\\
{51}&{-124}&{19}&{144}&{-71}&{-52}&{33}&0\\
{-20}&{50}&{-12}&{-52}&{28}&{18}&{-12}&0\\
5&{-12}&1&{16}&{-9}&{-4}&3&0\\
0&{-2}&8&{-12}&8&{-2}&0&0
\endmatrix\right),$$
$$V_2^{\sim\ast}(2)=8\left(\matrix
{400}&{-1243}&{972}&{542}&{-1028}&{357}&0&0\\
{-259}&{808}&{-642}&{-332}&{653}&{-228}&0&0\\
{156}&{-489}&{396}&{186}&{-384}&{135}&0&0\\
{-85}&{268}&{-222}&{-92}&{203}&{-72}&0&0\\
{40}&{-127}&{108}&{38}&{-92}&{33}&0&0\\
{-15}&{48}&{-42}&{-12}&{33}&{-12}&0&0\\
4&{-13}&{12}&{2}&{-8}&3&0&0\\
{-1}&4&{-6}&4&{-1}&0&0&0
\endmatrix\right),$$
$$V_2^{\sim\ast}(3)=8\left(\matrix
{119}&{-476}&{714}&{-476}&{119}&0&0&0\\
{-76}&{304}&{-456}&{304}&{-76}&0&0&0\\
{45}&{-180}&{270}&{-180}&{45}&0&0&0\\
{-24}&{96}&{-144}&{96}&{-24}&0&0&0\\
{11}&{-44}&{66}&{-44}&{11}&0&0&0\\
{-4}&{16}&{-24}&{16}&{-4}&0&0&0\\
1&{-4}&6&{-4}&1&0&0&0\\
0&0&0&0&0&0&0&0
\endmatrix\right).$$
The above matices $S^\sim_l$ and $V_l^{\sim\ast}(i)$
have the following properties:
\begin{equation}\label{eq:10}
A_l^\sim(z;-\nu)T_{4+2l,-1}A_{l,k}(z;\nu)=T_{4+2l,-1},
\end{equation}
\begin{equation}\label{eq:11}
S^\sim_lT_{4+2l,-1}=(S^\sim_lT_{4+2l,-1})^{-1}
\end{equation}
\begin{equation}\label{eq:12}
S^\sim_lT_{4+2l,-1}V_l^{\sim\ast}(i)=
-(-1)^iV_l^{\sim\ast}(i)T_{4+2l,-1}S^\sim_l,
\end{equation}
\begin{equation}\label{eq:13}
V_l^{\sim\ast}(i)T_{4+2l,-1}V_l^{\sim\ast}(k)=0T_{4+2l,-1},
\end{equation}
where
$$l=0,\,1,\,2,\,i\in[0,1+l]\cap{\mathbb Z},\,k\in[0,1+l]\cap{\mathbb Z}.$$ 
\bf Proof. \rm Full proof can be found in [2] -- [6].
{\begin{center}\large\bf References.\end{center}}
\footnotesize
\vskip4pt
\refstepcounter{r}\noindent[\ther]
R.Ap\'ery, Interpolation des fractions continues\\
\hspace*{3cm}  et irrationalite de certaines constantes,\\
\hspace*{3cm} Bulletin de la section des sciences du C.T.H., 1981, No 3,
 37 -- 53;
\label{r:cd}\\
\refstepcounter{r}\noindent [\ther]
 L.A.Gutnik, On some systems of difference
equations. Part 1.\\
\hspace*{3cm} Max-Plank-Institut f\"ur Mathematik, Bonn,\\
\hspace*{3cm} Preprint Series, 2006, 23, 1 -- 37.
\label{r:caj3}\\
\refstepcounter{r}
\noindent [\ther] \rule{2cm}{.3pt}, On some systems of difference
equations. Part 2.\\
\hspace*{3cm} Max-Plank-Institut f\"ur Mathematik, Bonn,\\
\hspace*{3cm} Preprint Series, 2006, 49, 1 -- 31.
\label{r:caj4}\\
\refstepcounter{r}
\noindent [\ther] \rule{2cm}{.3pt}, On some systems of difference
equations. Part 3.\\
\hspace*{3cm} Max-Plank-Institut f\"ur Mathematik, Bonn,\\
\hspace*{3cm} Preprint Series, 2006, 91, 1 -- 52.
\label{r:caj5}\\
\refstepcounter{r}
\noindent [\ther] \rule{2cm}{.3pt}, On some systems of difference
equations. Part 4.\\
\hspace*{3cm} Max-Plank-Institut f\"ur Mathematik, Bonn,\\
\hspace*{3cm} Preprint Series, 2006, 101, 1 -- 49.
\label{r:caj6}\\
\refstepcounter{r}
\noindent [\ther] \rule{2cm}{.3pt}, On some systems of difference
equations. Part 5.\\
\hspace*{3cm} Max-Plank-Institut f\"ur Mathematik, Bonn,\\
\hspace*{3cm} Preprint Series, 2006,115, 1 -- 9. 
\label{r:caj7}
\vskip 10pt
 {\it E-mail:}{\sl\ gutnik$@@$gutnik.mccme.ru}
\end{document}